# NUMERICAL METHOD FOR SPACE-TIME FRACTIONAL DIFFUSION: A STOCHASTIC APPROACH[*]

TENGTENG CUI[†], CHENGTAO SHENG[‡], BIHAO SU[§], AND ZHI ZHOU[†]

**Abstract.** In this paper, we develop and analyze a stochastic algorithm for solving space-time fractional diffusion models, which are widely used to describe anomalous diffusion dynamics. These models pose substantial numerical challenges due to the memory effect of the time-fractional derivative and the nonlocal nature of the spatial fractional Laplacian and the, leading to significant computational costs and storage demands, particularly in high-dimensional settings. To overcome these difficulties, we propose a Monte Carlo method based on the Feynman–Kac formula for space-time fractional models. The novel algorithm combines the simulation of the monotone path of a stable subordinator in time with the "walk-on-spheres" method that efficiently simulates the stable Lévy jumping process in space. We rigorously derive error bounds for the proposed scheme, explicitly expressed in terms of the number of simulation paths and the time step size. Numerical experiments confirm the theoretical error bounds and demonstrate the computational efficiency of the method, particularly in domains with complex geometries or high-dimensional spaces. Furthermore, both theoretical and numerical results emphasize the robustness of the proposed approach across a range of fractional orders, particularly for small fractional values, a capability often absent in traditional numerical methods.

**1. Introduction.** In this work, we focus on the numerical treatment of the following initial/boundary value problem for the space-time fractional diffusion equation in a bounded domain $\Omega \subseteq \mathbb{R}^n$, where $n \geq 2$, $\beta \in (0, 1]$ and $\alpha \in (0, 2]$:

$$(1.1) \quad \begin{cases} \partial_t^\beta u(t, x) + (-\Delta)^{\frac{\alpha}{2}} u(t, x) = f(t, x), & (t, x) \in (0, T] \times \Omega, \\ u(t, x) = g(t, x), & (t, x) \in [0, T] \times \Omega^c, \\ u(t, x) = u_0(x), & x \in \Omega. \end{cases}$$

Here, $\partial_t^\beta u$ denotes the Djrbashian–Caputo fractional derivative of order $\beta \in (0, 1)$ in time, which is defined by [22, Definition 2.3]:

$$(1.2) \quad \partial_t^\beta u(t) = \frac{1}{\Gamma(1 - \beta)} \int_0^t (t - s)^{-\beta} u'(s) \, \mathrm{d}s,$$

where $\Gamma(z)$ is the Euler gamma function. When $\beta \to 1^-$, the fractional derivative $\partial_t^\beta u$ reduces to the first-order derivative in time, $u'(t)$. The term $(-\Delta)^{\frac{\alpha}{2}} u$ represents the fractional Laplacian, which has the following hypersingular integral representation [10]:

$$(1.3) \quad (-\Delta)^{\frac{\alpha}{2}} u(x) = c_{n,\alpha} \text{ p.v.} \int_{\mathbb{R}^n} \frac{u(x) - u(y)}{|x - y|^{n+\alpha}} \, \mathrm{d}y \quad \text{with} \quad c_{n,\alpha} := \frac{2^\alpha \Gamma\left(\frac{n+\alpha}{2}\right)}{\pi^{\frac{n}{2}} \Gamma\left(1 - \frac{\alpha}{2}\right)},$$

---

[*]The work of C. Sheng is supported by the National Natural Science Foundation of China (Projects 12201385 and 12271365). The work of B. Su is supported by the Scientific Research Start-up Funds of Hainan University (Projects XJ2400010491). The work of Z. Zhou is supported by by National Natural Science Foundation of China (Project 12422117), Hong Kong Research Grants Council (15303122) and an internal grant of Hong Kong Polytechnic University (Project ID: P0038888, Work Programme: 1-ZVX3).

[†]Department of Applied Mathematics, The Hong Kong Polytechnic University, Kowloon, Hong Kong (tengcui@polyu.edu.hk, zhizhou@polyu.edu.hk).

[‡]School of Mathematics, Shanghai University of Finance and Economics, Shanghai 200433, China(ctsheng@sufe.edu.cn)

[§]School of Mathematics and Statistics, Hainan University, Haikou 570100, China (bihaosu@hainanu.edu.cn)

where "p.v." indicates the principal value of the integral, and $c_{n,\alpha}$ is a normalization constant. This pointwise expression can also be derived from the following Fourier transform:

$$(-\Delta)^{\frac{\alpha}{2}} u(x) = \mathscr{F}^{-1}\big[|\xi|^\alpha \mathscr{F}[u](\xi)\big](x), \quad x \in \mathbb{R}^n,$$

where $\mathscr{F}$ and $\mathscr{F}^{-1}$ denote the Fourier and inverse Fourier transforms, respectively. In the special case of $\alpha = 2$, the fractional Laplacian simplifies to the standard negative Laplacian $-\Delta$, and the volume constraint boundary condition $u = g$ on $\Omega^c$ reduces to the classical Dirichlet boundary condition $u = g$ on $\partial\Omega$.

The space-time fractional evolution models (1.1) are widely used to describe anomalous diffusion, where the mean squared displacement of particles grows either slower (subdiffusion) or faster (superdiffusion) than in normal diffusion. These models have garnered significant attention over the past few decades due to their broad range of applications, including protein diffusion within cells [18], contaminant transport in groundwater [26], thermal diffusion in media with fractal geometry [35], and heat conduction with memory effects [47], among others. These fractional operators are closely linked to stochastic processes that generalize classical Brownian motion, providing a probabilistic interpretation of the underlying dynamics. For instance, the fractional Laplacian $(-\Delta)^{\frac{\alpha}{2}}$ is the generator of symmetric $\alpha$-stable Lévy processes, which model random jumps or discontinuous paths rather than continuous trajectories, making them suitable for describing superdiffusive behavior [34, 46]. Similarly, the Djrbashian–Caputo fractional derivative $\partial_t^\beta$ is associated with subordinated stochastic processes, where the time evolution is governed by an inverse stable subordinator, capturing the memory effects that characterize subdiffusive behavior [33, 8, 13]. It has been extensively reported in the literatures that, the inclusion of fractional operators has a profound impact on the diffusive dynamics [2, 14, 30, 44, 12] and solution regularity [37, 3, 36].

Over the past two decades, significant progress has been achieved in the development and analysis of numerical methods for solving space-time nonlocal and fractional models. For a comprehensive and up-to-date overview, refer to the monograph [24], the survey articles [6, 9, 23, 32], and the references therein. One of the main challenges in these models arises from the spatial nonlocality and memory effects introduced by the fractional operators, which lead to high computational costs and substantial storage demands. To address these difficulties, many fast algorithms have been developed. Notable examples include fast and oblivious convolution quadrature methods [38, 4, 16] and sum-of-exponentials approximations [5, 20, 50] for time-fractional evolution models, as well as preconditioned solvers based on Toeplitz-like systems [31, 11] for fractional-in-space models. See also [49, 48, 29] for some parallel-in-time implementations of numerical schemes for solving fractional models. Another significant challenge arises from the stronger singularity behavior at smaller fractional orders $\alpha$ and $\beta$. In such cases, it is often necessary to refine the discretization meshes, particularly near the initial time and domain boundaries, by employing highly graded meshes to achieve optimal convergence rates (see, e.g., [1, 15, 7, 42]). However, this refinement can significantly amplify round-off errors, resulting in unstable computations and further complicating the numerical treatment of these problems (see, e.g., [41, Remark 11]).

Unlike the aforementioned deterministic numerical methods, stochastic algorithms apply the Feynman–Kac formula, which establishes a connection between PDE models and stochastic processes. This connection enables the use of Monte Carlo methods, making these algorithms particularly well-suited for solving problems involving nonlocal operators and high-dimensional settings. In recent years, some effort has been invested in exploring stochastic simulations for nonlocal problems. EKolokoltsov et al. [27] developed an effi-cient Monte Carlo method for the time-fractional diffusion models, i.e. (1.1) with $\alpha = 2$, based on the exact probabilistic representation, and provided *a priori* error estimates. See also the probabilistic representation for general nonlocal-in-time diffusion models in [8, 13]. For space-fractional models, Kyprianou et al. [28] proposed a walk-on-spheres algorithm for the integral fractional Laplacian in two dimensions, and then Shardlow extended this to fractional eigenvalue problems in two dimensions via multilevel Monte Carlo methods [39]. While these methods can be applied to high-dimensional settings, the complexity of simulating $\alpha$-stable processes in higher dimensions poses significant challenges. Jiao et al. [21] modified this strategy by proposing a quadrature rule to evaluate integral representation in the ball and applying rejection sampling method, to extend the method to a general domain. To further simplify the high-dimensional quadrature, Sheng et al. [40] applied spherical coordinates and the explicit expressions of the Green's function and the Poisson kernel to derive the jump distances. This innovation significantly reduced computational complexity, and their work provided rigorous error estimates for the stochastic algorithm. Despite these advances, no Monte Carlo-based stochastic algorithm has yet been developed to solve space-time fractional diffusion equations (1.1), likely due to the intricate interaction between the temporal evolution driven by the inverse stable subordinator and the spatial evolution governed by the Lévy jumping process.

In this work, we develop an easy-to-implement stochastic algorithm for the space-time fractional diffusion model. The core idea of the algorithm is to apply the Monte Carlo approximation based on the associated Feynman–Kac formula [45] and to effectively combine two key components: the simulation of the monotone path of a stable subordinator in time, corresponding to the time-fractional derivative[27], and the "walk-on-spheres" method, which efficiently simulates the stable Lévy jumping process in space associated with the fractional Laplacian [40]. Compared to classical discretization schemes, the proposed method offers several notable advantages:

(1) **Efficient simulation**: The algorithm simulates a jumping process over time and efficiently tracks particle trajectories in space. Unlike classical numerical methods, such as finite element or finite difference approaches, which struggle with the computational challenges of nonlocality, this algorithm significantly reduces computational and storage costs.

(2) **Unified framework**: The algorithm acts as a unified stochastic solver for both normal and anomalous diffusion models. It naturally simplifies to stochastic solvers for subdiffusion models (as $\alpha \to 2^-$), superdiffusion models (as $\beta \to 1^-$), and the normal diffusion model (when $\alpha \to 2^-$ and $\beta \to 1^-$ simultaneously).

(3) **Robustness for small fractional-orders**: The algorithm demonstrates robustness for small fractional orders $\alpha$ and $\beta$, which are particularly challenging for classical methods. For example, Figures 4-6 in Section 4 illustrate the accuracy and robustness of the proposed methods even for very small fractional orders $\alpha$ and $\beta$.

(4) **Meshless framework**: The proposed scheme is spatially meshless, making it well-suited for irregular domains with complex geometries. For example, Example 4 in Section 4 illustrates a simulation performed on a hexagonal hailstone-shaped domain.

(5) **Scalability for high dimensions**: The Monte Carlo approximation is inherently parallelizable for each trajectory, making it particularly well-suited for high-dimensional simulations. Furthermore, its convergence rate is independent of the spatial dimension, ensuring robust performance without suffering from the curse of dimensionality. For example, Figure 4 illustrates the approximation errors in a 100-dimensional problem, demonstrating its efficiency and practicality in such a challenging setting.

Moreover, we establish an error estimate for the proposed stochastic scheme. Specifically, under some mild smoothness condition on problem data, and let $u$ and $u_M^*$ denote the solutions of (1.1) and the numerical scheme (3.3), respectively. Then, the following error estimate holds (cf. Theorem 3.4)

$$\left(\mathbb{E}\left[\|u(t,\cdot) - u_M^*(t,\cdot)\|_{L^2(\Omega)}^2\right]\right)^{\frac{1}{2}} \leq C\left((\Delta t)^{(1-\epsilon)/2} + M^{-1/2}\right),$$

for any arbitrarily small $\epsilon > 0$. Here $M$ is the sampling number in the Monte Carlo approximation, and $\Delta t$ is the step size used in the discretization for the time integration of the source term. It is important to note that the constant $C$ is independent of fractional orders $\beta \in (0,1]$ and $\alpha \in (0,2]$, which indicates the robustness of the proposed scheme for the extreme cases when $\alpha, \beta$ are close to zero.

The rest of the paper is organized as follows. In section 2, we introduce some notations and provides some necessary preliminaries for Feynman–Kac formulas. In Section 3, we construct the full-discrete Monte Carlo scheme for solving the space-time fractional model (1.1) and provide a rigorous error analysis. We give in Section 4 several numerical examples to verify the theoretical results and demonstrate the efficiency of the proposed scheme. Some concluding remarks are given in the last section.

## 2. Preliminaries.
In this section, we introduce some notation and review some basic results about several stochastic processes closely related to the fractional differential operators which will serve as a fundamental component in the construction of the stochastic algorithm.

### 2.1. Stochastic processes and Feynman–Kac formula.
We begin by introducing the symmetric $\alpha$-stable Lévy process, a stochastic process characterized by heavy-tailed distributions and jumps, with the fractional Laplacian operator serving as its infinitesimal generator. For $\alpha \in (0,2)$, we denote by $X^\alpha$ the rotationally symmetric $\alpha$-stable Lévy process in $\mathbb{R}^n$. In the special case where $\alpha = 2$, the process $X^\alpha$ reduces to $\sqrt{2}B_s$, where $B_s$ represents the standard Brownian motion. A key property of the symmetric $\alpha$-stable Lévy process is its $1/\alpha$-self-similarity: for any $k > 0$, the processes $\{X^\alpha(ks); s \geq 0\}$ and $\{k^{1/\alpha}X^\alpha(s); s \geq 0\}$ share the same finite-dimensional distributions. This self-similarity, combined with its jump-driven dynamics, makes the $\alpha$-stable Lévy process a natural generalization of Brownian motion for modeling heavy-tailed phenomena.

At any fixed time $t$, the process $X^\alpha(t)$ is an $\alpha$-stable random variable, whose probability distribution is fully characterized by its characteristic function. In general, the characteristic function of an $\alpha$-stable random variable is [19]

$$\ln \psi(\theta) = \begin{cases} -\sigma^\alpha |\theta|^\alpha \left[1 - i\rho \, \text{sgn}(\theta) \tan\left(\frac{\alpha\pi}{2}\right)\right] + i\mu\theta, & \text{if } \alpha \in (0,1) \cup (1,2], \\ -\sigma|\theta|\left[1 + i\rho\frac{2}{\pi} \text{sgn}(\theta) \ln|\theta|\right] + i\mu\theta, & \text{if } \alpha = 1, \end{cases}$$

where $\alpha \in (0,2]$ is the stability index, $\rho \in [-1,1]$ is the skewness parameter, $\mu \in \mathbb{R}$ is the shift parameter, and $\sigma > 0$ is the scale parameter. For special cases, when $\rho = 0$, the random variable is symmetric. For example, the case $\alpha = 1$ corresponds to the Cauchy distribution, while the case $\alpha = 2$ corresponds to the Gaussian distribution $\mathcal{N}(\mu, \Sigma)$. Specifically, the characteristic function of a symmetric $\alpha$-stable random variable corresponding to fractional Laplacian (1.3) is:

$$(2.1) \qquad \psi(\theta) = \mathbb{E}[\exp(\mathrm{i}\theta \cdot X^\alpha(t))] = \exp(-t|\theta|^\alpha) \quad \text{with } \alpha \in (0,2], \ \theta \in \mathbb{R}^n.$$

A subordinator is a special type of stochastic process that is an increasing Lévy process. A Lévy process on $\mathbb{R}$ is said to be increasing if $Y^\beta(s)$ increases as a function of $s$. When a Lévy process has this monotonic property, it is referred to as a subordinator. In particular, $Y^\beta$ is called a $\beta$-stable subordinator for $\beta \in (0,1)$, where $\beta$ is the stability index. A $\beta$-stable subordinator starts at zero, has zero drift, and satisfies the scaling property $Y^\beta(s) = s^{1/\beta} Y^\beta(1)$ for all $s > 0$, where $Y^\beta(1)$ follows the stable distribution $S_\beta(1,1,0)$. The Laplace exponent of the process is given by:

$$(2.2) \qquad \mathbb{E}[e^{-kY^\beta(s)}] = e^{-sk^\beta}, \quad k, s > 0,$$

which characterizes the process as having a stable distribution with index $\beta$.

Next, we introduce the Feynman–Kac formula associated with space-time fractional equations on bounded domains. These formulas serve as the foundation and a pivotal starting point for the development of our stochastic algorithms. There exists a unique stochastic representation for the problem (1.1):

$$
\begin{aligned}
(2.3) \quad u(t,x) = {}& \mathbb{E}_{X^\alpha(0)=x} \left[ u_0 \left( X^\alpha(\tau_t) \right) \mathbb{I}_{\{\tau_t < \tau_\Omega(x)\}} + g \left( t - Y^\beta(\tau_\Omega(x)), X^\alpha(\tau_\Omega(x)) \right) \mathbb{I}_{\{\tau_t \geq \tau_\Omega(x)\}} \right] \\
& + \mathbb{E}_{X^\alpha(0)=x} \left[ \int_0^{\tau_t \wedge \tau_\Omega(x)} f \left( t - Y^\beta(s), X^\alpha(s) \right) ds \right],
\end{aligned}
$$

where $Y^\beta$ is the standard $\beta$-stable subordinator, and $X^\alpha$ is the symmetric $\alpha$-stable Lévy process, as defined in (2.2) and (2.1), respectively. Here,

$$(2.4) \qquad \tau_\Omega(x) = \inf\{s > 0 \mid X^\alpha(s) \notin \Omega \text{ with } X^\alpha\}$$

is the *first exit time* of $X^\alpha$ from the domain $\Omega$, and

$$(2.5) \qquad \tau_t = \inf\{s > 0 \mid t - Y^\beta(s) \leq 0, \ t \in (0,T]\}$$

is the *first exit time* where $t - Y^\beta(s)$ across origin.

The formula (2.3) establishes a direct connection between the solution and stochastic processes. These insights motivate us to harness the power of stochastic representations to develop efficient stochastic algorithms for solving space-time fractional equations on bounded domains.

**2.2. Numerical simulation of the motion paths in space/time.** To simulate $\alpha$-stable processes in high-dimensional space, we employ the walk-on-spheres method, which was proposed recently in [40, 43] for fractional Laplacian. For a given time partition, the trajectory of the process $X^\alpha$ is approximated using a sequence of spheres, where the radius of each sphere is determined by the time step size of the partition.

To this end, let $c_0 \in \mathbb{R}^n$ and define $r > 0$ as the radius of a ball inscribed in $\mathbb{R}^n$, centered at $c_0$. We denote this ball as $\mathbb{B}_r^n(c_0) = \{x \in \mathbb{R}^n : |x - c_0| \leq r\}$. For simplicity, let $\mathbb{B}_r^n = \mathbb{B}_r^n(0)$. We then define the exit time of the process $X^\alpha$ from $\mathbb{B}_r^n$ as

$$(2.6) \qquad \tau_r = \inf\{s : X^\alpha(s) \notin \mathbb{B}_r^n \mid X^\alpha(0) = 0\}.$$

The expectation of the first exit time of the $\alpha$-stable Lévy process, which is very useful in algorithm development, has the following explicit expression [17, eq. (A) and eq. (B)]:

$$(2.7) \qquad \mathbb{E}[\tau_r] = r^\alpha C_n^\alpha, \quad \mathbb{E}[(\tau_r)^2] = \alpha r^\alpha (C_n^\alpha)^2 \int_0^{r^2} \nu^{\frac{\alpha}{2}-1} {}_2F_1(-\frac{\alpha}{2}, \frac{n}{2}; \frac{n+\alpha}{2}; \nu r^{-2}) d\nu,$$

where the constant

$$(2.8) \qquad C_n^\alpha = \frac{\Gamma(\frac{n}{2})}{2^\alpha \Gamma(1 + \frac{\alpha}{2})\Gamma(\frac{n+\alpha}{2})}.$$

and the hypergeometric function takes the form as

$$(2.9) \qquad {}_2F_1(a, b; c, z) = \frac{\Gamma(c)}{\Gamma(b)\Gamma(c - b)} \int_0^1 t^{b-1}(1 - t)^{c-b-1}(1 - zt)^{-a} \mathrm{d}t.$$

Then for any $\epsilon > 0$, there holds [43, eq. (13)]

$$(2.10) \qquad \lim_{r \to 0} \mathbb{P}(|\tau_r - \mathbb{E}[\tau_r]| \leq \epsilon) > 1 - \frac{r^{2\alpha}(\tilde{C}_n^\alpha)^2}{\epsilon^2} \to 1^- \quad \text{as} \quad r \to 0.$$

REMARK 2.1. *It is straightforward to observe that*

$$\begin{aligned}
{}_2F_1(-\frac{\alpha}{2}, \frac{n}{2}; \frac{n+\alpha}{2}; \nu r^{-2}) &= \frac{\Gamma(\frac{n+\alpha}{2})}{\Gamma(\frac{n}{2})\Gamma(\frac{\alpha}{2})} \int_0^1 t^{\frac{n}{2}-1}(1-t)^{\frac{\alpha}{2}-1}(1 - zt)^{\frac{\alpha}{2}} \mathrm{d}t \\
&\leq \frac{\Gamma(\frac{n+\alpha}{2})}{\Gamma(\frac{n}{2})\Gamma(\frac{\alpha}{2})} \int_0^1 t^{\frac{n}{2}-1}(1-t)^{\frac{\alpha}{2}-1} \mathrm{d}t = 1
\end{aligned}$$

*As a result, for $n \geq 2$ and $\alpha \in (0, 2]$, the random variable $\tau_r$ has a finite second moment:*

$$\mathbb{E}[(\tau_r)^2] \leq \alpha r^\alpha (C_n^\alpha)^2 \int_0^{r^2} \nu^{\frac{\alpha}{2}-1} \mathrm{d}\nu \leq 2^{1-2\alpha} r^{2\alpha} \left(\frac{\Gamma(\frac{n}{2})}{\Gamma(1 + \frac{\alpha}{2})\Gamma(\frac{n+\alpha}{2})}\right)^2 \leq 4r^{2\alpha}.$$

In practice, to simulate the $\beta$-stable subordinator in time, we fix a small time step size $\Delta t$ and partition the time interval $I := [0, \infty)$ into subintervals with the grid

$$(2.11) \qquad \mathcal{T} := \{0 = t_0 < t_1 < t_2 < \cdots\}, \quad \text{where} \quad t_i = i\Delta t.$$

Next, we generate a path of $Y^\beta(t)$ at the grid points using the recurrence relation:

$$(2.12) \qquad Y^\beta(t_i) \sim Y^\beta(t_{i-1}) + (\Delta t)^{1/\beta} \eta_i, \quad \text{where} \quad \eta_i \sim S_\beta(1, 1, 0).$$

We then define the stopping time as

$$\tau_t = \Delta t \times \min\{i : Y^\beta(t_i) \geq t\}.$$

It is worth noting that, theoretically, for a given $t \in (0, T]$, the stopping time satisfies the following relationship (see [27, Lemma 2.1]):

$$\tau_t \stackrel{d}{=} (t/\eta)^\beta, \quad \text{where} \quad \eta \sim S_\beta(1, 1, 0).$$

From the relation $\mathbb{E}[\tau_r] = \Delta t$, the first relation in (2.7) implies that the radius of the ball is given by:

$$(2.13) \qquad r = (\Delta t / C_n^\alpha)^{1/\alpha}.$$

Then we can approximate the trajectory of the process $X^\alpha(t)$ on a time grid by simulating it as a sequence of small balls, where the radius $r$ is determined by the time step size $\Delta t$. Specifically, in each time interval $[t_{i-1}, t_i)$, the process $X^\alpha(t)$ starts at the location $X^\alpha(t_{i-1})$. For $t \in [t_{i-1}, t_i)$, the trajectory of $X^\alpha(t)$ remains confined within the ball $\mathbb{B}_r^n(X^\alpha(t_{i-1}))$. At time $t_i$, the stochastic process $X^\alpha(t)$ undergoes a jump and leaves the ball $\mathbb{B}_r^n(X^\alpha(t_{i-1}))$, with $\tau_r \approx \Delta t$ according to (2.10). Hence, the motion path of the process $X^\alpha(t)$ can be effectively approximated by the motion path of the sequence of balls; see Figure 1.

Since the radius of the ball is fixed, the primary focus will be on simulating the position of the center of each ball. Directly working with the relevant quantities in a higher-dimensional Cartesian coordinate system is a challenging task. Therefore, for simplicity and convenience in calculations, we propose performing the computations in the spherical coordinate system in dimension $n$:

$$(2.14) \quad \begin{aligned} & x_1 = \rho \cos\theta_1; \ x_2 = \rho \sin\theta_1 \cos\theta_2; \ \cdots\cdots; \ x_{n-1} = \rho \sin\theta_1 \cdots \sin\theta_{n-2} \cos\theta_{n-1}; \\ & x_n = \rho \sin\theta_1 \cdots \sin\theta_{n-2} \sin\theta_{n-1}, \ \text{with } \theta_1, \cdots, \theta_{n-2} \in [0, \pi], \ \text{and } \theta_{n-1} \in [0, 2\pi], \end{aligned}$$

In the spherical coordinate system, the angular components are uniformly distributed across all directions, so we only need to determine the radial jump length.

For a given time partition $\mathcal{T}$ in (2.11) and the radius $r$ in (2.13), one can derive the following identity on the position of $X^\alpha(t_i)$

$$(2.15) \quad X^\alpha(t_i) = X^\alpha(t_{i-1}) + J_i \cdot \begin{bmatrix} \cos\theta_1 \\ \sin\theta_1 \cos\theta_2 \\ \cdots\cdots \\ \sin\theta_1 \cdots \sin\theta_{n-2} \sin\theta_{n-1} \end{bmatrix},$$

where the $i$-th jump distance is given by [40, Lemma 2.3]

$$(2.16) \quad J_i = \sqrt{\dfrac{r^2}{B\left(1 - \frac{\alpha}{2}, \frac{\alpha}{2}\right) - B^{-1}\left(\frac{\pi\omega}{\sin(\pi\alpha/2)}; 1 - \frac{\alpha}{2}, \frac{\alpha}{2}\right)}}, \quad \omega \in (0,1),$$

and $B^{-1}(\cdot; a, b)$ denotes the inverse function of incomplete beta function $B(\cdot; a, b)$, $B(a, b) = B(1; a, b)$ denotes the beta function.

In practical calculations, we only need to sample $\omega$ as a uniformly distributed random number over $(0, 1)$ and substitute it into (2.16) to obtain the jump distance $J_i$. Similarly, we sample a uniformly distributed random number for the spherical direction and combine it with (2.14) to determine the value of the position $X^\alpha(\tau_i)$. Hence, the motion path of the stochastic process can be approximated by the motion path of the balls, as shown in Figure 1.

**3. The Monte Carlo scheme and its error estimate.** In this section, we describe the implementations of the Monte Carlo method for the solution of (1.1) using the Feynman–Kac formula (2.3), and provide rigorous error estimate of the proposed scheme.

**3.1. The implementation of Monte Carlo method.** For any fixed $t$, we generate the path of $Y^\beta$ using the relation (2.12) as well as the exit time $\tau_t$. Let $N = \tau_t / \Delta t$. We also generate a sequence of small balls with a radius given by $r = (\Delta t / C_n^\alpha)^{1/\alpha}$. For any given $x \in \Omega$, we set $X^\alpha(0) = x$ and draw a ball $\mathbb{B}_r^n(X^\alpha(0))$ within the domain $\Omega$. Using (2.15), we generate a sequence of balls centered at $X^\alpha(t_i)$ for $i = 1, 2, \ldots, N$. The remaining task is to determine whether the points $X^\alpha(t_i)$ lie inside the domain.

(i) If $X^\alpha(t_i) \in \Omega$ for all $i = 1, 2, \ldots, N$, then we conclude that $\tau_\Omega > \tau_t$. Then we apply the right-point quadrature rule to approximate the integral of (2.3). Specifically,

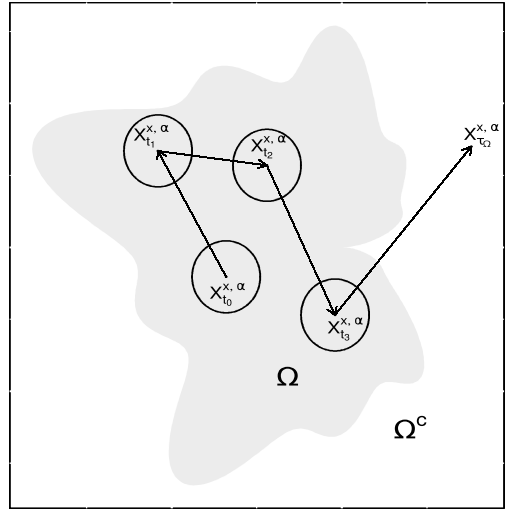

Fig. 1. *The motion trajectory of $\alpha$-stable process on $\Omega$*

we define

$$(3.1) \qquad S_j\left(\tau_t; x\right) = u_0\left(X_j^{\alpha}(\tau_t)\right) + \Delta t \sum_{i=1}^{N_j} f(t - Y_j^{\beta}(t_i), X_j^{\alpha}(t_i)).$$

where the index $j$ represents the $j$-th experiment.

(ii) Otherwise, if there exists $k$ such that $1 \leq k \leq N$ and $X^{\alpha}(t_k) \notin \Omega$, then we conclude that $\tau_{\Omega} = t_{k_*} \leq \tau_t$ with $k_* = \min\{k : X^{\alpha}(t_k) \notin \Omega\}$, and obtain the following approximation

$$(3.2) \qquad S_j\left(\tau_{\Omega}, x\right) = g\left(t - Y_j^{\beta}(\tau_{\Omega}), X_j^{\alpha}(\tau_{\Omega})\right) + \Delta t \sum_{i=1}^{k_*} f(t - Y_j^{\beta}(t_i), X_j^{\alpha}(t_i)).$$

A combination of the above two equations and the law of large numbers, we can obtain the numerical solution by using the following full discretization scheme

$$
\begin{aligned}
(3.3) \qquad u_M^*(t, x) &= \frac{1}{M} \sum_{j=1}^{M} S_j\left(\tau_t \wedge \tau_{\Omega}, x\right) \\
&= \frac{1}{M} \sum_{j=1}^{M} \Big[ u_0\left(X_j^{\alpha}(\tau_t)\right) \mathbb{I}_{\{\tau_t < \tau_{\Omega}(x)\}} + g\left(t - Y_j^{\beta}(\tau_{\Omega}), X_j^{\alpha}(\tau_{\Omega})\right) \mathbb{I}_{\{\tau_t \geq \tau_{\Omega}(x)\}} \\
&\qquad + \Delta t \sum_{i=1}^{k_j} f(t - Y_j^{\beta}(t_i), X_j^{\alpha}(t_i)) \Big].
\end{aligned}
$$

where $M$ represents the total number of experiments and

$$(3.4) \qquad k_j = \begin{cases} N_j, & \text{if } \tau_t < \tau_{\Omega}(x); \\ k_{*,j} = \min\{k : X_j^{\alpha}(t_k) \notin \Omega\}, & \text{if } \tau_{\Omega}(x) < \tau_t. \end{cases}$$

**3.2. Error estimate.** Next, we aim to derive an error estimate for the Monte-Carlo scheme (3.3). First, we assume that problem data $u_0$, $g$, $f$ uniformly bounded, i.e.,

$$\text{(3.5)} \quad \|u_0\|_{L^\infty(\Omega)} \leq K_0, \ \|g(t,\cdot)\|_{L^\infty(\Omega^c)} \leq K_g \text{ and } \|f(t,\cdot)\|_{L^\infty(\Omega)} \leq K_f, \quad \text{for all } t \in (0,T].$$

Moreover, we assume that $f$ satisfies, for some small $\epsilon > 0$

$$\text{(3.6)} \quad \begin{aligned} |f(t,x) - f(s,x)| &\leq c_1 |t-s|^{\gamma_1} \text{ for all } t,s \in [0,T] \quad \text{with } \gamma_1 = \beta(1-\epsilon)/2; \\ |f(t,x) - f(t,y)| &\leq c_2 \|x-y\|^{\gamma_2} \text{ for all } x,y \in \overline{\Omega} \quad \text{with } \gamma_2 = \alpha(1-\epsilon)/2, \end{aligned}$$

Here we recall that the spatial symmetric $\alpha$-stable Lévy process process $X^\alpha$ satisfies the following Hölder property [25, Example 25.10]

$$\text{(3.7)} \quad \mathbb{E}[\|X^\alpha(s_1) - X^\alpha(s_2)\|^p] \leq \frac{2^p \Gamma(\frac{1+p}{2})\Gamma(1-\frac{p}{\alpha})}{\sqrt{\pi}\Gamma(1-\frac{p}{2})}(s_1-s_2)^{p/\alpha}, \quad \text{for any } p \in (0,\alpha).$$

Similarly, the temporal $\beta$-stable subordination process $Y^\beta$ satisfies

$$\text{(3.8)} \quad \mathbb{E}[|Y^\beta(t)|^p] = \frac{2^p \Gamma(\frac{1+p}{2})\Gamma(1-\frac{p}{\beta})}{\sqrt{\pi}\Gamma(1-\frac{p}{2})}t^{p/\beta} \quad \text{for any } p \in (0,\beta).$$

Let $u$ and $u_M^*$ be the solution of (2.3) and (3.3), respectively. Then, we have the following splitting

$$\text{(3.9)} \quad \|(u - u_M^*)(t)\|_{L^2(\Omega)} \leq \sum_{i=1}^3 B_i,$$

where $B_i$ are defined by

$$\begin{aligned} B_1 =& \Big\| \mathbb{E}_{X^\alpha(0)=x}\Big[u_0\left(X^\alpha(\tau_t)\right)\mathbb{I}_{\tau_\Omega(x) > \tau_t} + g\left(t - Y^\beta(\tau_\Omega), X^\alpha(\tau_\Omega)\right)\mathbb{I}_{\tau_\Omega(x) \leq \tau_t}\Big] \\ & - \frac{1}{M}\sum_{j=1}^M \Big[u_0\left(X_j^\alpha(\tau_{t,j})\right)\mathbb{I}_{\tau_{\Omega,j}(x) > \tau_{t,j}} + g\left(t - Y_j^\beta(\tau_{\Omega,j}), X_j^\alpha(\tau_{\Omega,j})\right)\mathbb{I}_{\tau_{\Omega,j}(x) \leq \tau_{t,j}}\Big]\Big\|_{L^2(\Omega)}, \\ B_2 =& \Big\| \mathbb{E}_{X^\alpha(0)=x}\Big[\int_0^{\tau_t \wedge \tau_\Omega(x)} f(t - Y^\beta(s), X^\alpha(s))\mathrm{d}s\Big] \\ & - \frac{1}{M}\sum_{j=1}^M \int_0^{\tau_{t,j} \wedge \tau_{\Omega,j}(x)} f(t - Y_j^\beta(s), X_j^\alpha(s))\mathrm{d}s\Big\|_{L^2(\Omega)}, \\ B_3 =& \Big\| \frac{1}{M}\sum_{j=1}^M \int_0^{\tau_{t,j} \wedge \tau_{\Omega,j}(x)} f(t - Y_j^\beta(s), X_j^\alpha(s))\mathrm{d}s - \frac{1}{M}\sum_{j=1}^M \Big[\Delta t \sum_{i=1}^{k_j} f(t - Y_j^\beta(t_i), X_j^\alpha(t_i))\Big]\Big\|_{L^2(\Omega)}. \end{aligned}$$

where $k_j$ is the stopping index of $j$-th trajectory defined by

$$k_j = \begin{cases} \lfloor \tau_{t,j}/\Delta t \rfloor, & \text{if } \tau_{t,j} < \tau_{\Omega,j}(x); \\ \min\{k : X_j^\alpha(t_k) \notin \Omega\}, & \text{if } \tau_{\Omega,j}(x) < \tau_{t,j}. \end{cases}$$

Next, we estimate the bounds of the three terms $\{B_i\}_{i=1}^3$ in (3.9) one by one as follows.

LEMMA 3.1. *Suppose that the condition* (3.5) *holds valid. Then the term $B_1$ in* (3.9) *satisfies*

$$\left(\mathbb{E}[(B_1)^2]\right)^{1/2} \leq CM^{-1/2},$$

*where the constant*

$$C = \sqrt{2|\Omega| \max\left((K_0)^2, (K_g)^2\right)}.$$

*Note that the constant $C$ depends on $K_0$, $K_g$ (the bounds of the problem data) and $\Omega$ (the domain), but it does not depend on the parameters $\alpha \in (0, 2]$ and $\beta \in (0, 1]$, the time step $\Delta t$, or the number of samples $M$.*

*Proof.* Define the random variable

$$(3.10) \qquad Z = u_0\left(X^\alpha(\tau_t)\right) \mathbb{I}_{\tau_\Omega > \tau_t} + g\left(t - Y^\beta(\tau_\Omega), X^\alpha(\tau_\Omega)\right) \mathbb{I}_{\tau_\Omega \le \tau_t}.$$

Since $|u_0| \le L$ and $|g| \le K$, we have $|Z| \le \max(L, K) \triangleq C_b$. The Monte Carlo estimator for $\theta = \mathbb{E}[Z]$ is $\hat{\theta}_M = \frac{1}{M} \sum_{j=1}^{M} Z_j$, where $Z_j$ are i.i.d. copies of $Z$. The mean squared error of this estimator is given by

$$(3.11) \qquad \mathbb{E}\left[|\theta - \hat{\theta}_M|^2\right] = \frac{\mathrm{Var}(Z)}{M}.$$

Thus,

$$\mathbb{E}[(B_1)^2] \le |\Omega| \sup_{x \in \Omega} \mathbb{E}\left[|\theta - \hat{\theta}_M|^2\right] = |\Omega| \sup_{x \in \Omega} \frac{\mathrm{Var}(Z)}{M}$$

Now, using the facts that $\mathrm{Var}(Z) \le \mathbb{E}[Z^2]$ and $(a + b)^2 \le 2a^2 + 2b^2$, we obtain

$$\mathbb{E}[Z^2] \le 2\mathbb{E}\left[|u_0\left(X^\alpha(\tau_t)\right)|^2 \mathbb{I}_{\tau_\Omega > \tau_t}\right] + 2\mathbb{E}\left[|g\left(t - Y^\beta(\tau_\Omega), X^\alpha(\tau_\Omega)\right)|^2 \mathbb{I}_{\tau_\Omega \le \tau_t}\right]$$
$$\le 2(K_0)^2 \mathbb{P}(\tau_\Omega > \tau_t) + 2(K_g)^2 \mathbb{P}(\tau_\Omega \le \tau_t) \le 2 \max\left((K_0)^2, (K_g)^2\right),$$

where in the last inequality we apply the relation $\mathbb{P}(\tau_\Omega > \tau_t) + \mathbb{P}(\tau_\Omega \le \tau_t) = 1$. Let $C_Z = 2\max((K_0)^2, (K_g)^2)$. Then $\mathrm{Var}(Z) \le C_Z$, and

$$\left(\mathbb{E}[(B_1)^2]\right)^{1/2} \le CM^{-1/2} \text{ with } C = \sqrt{|\Omega|C_Z} = \sqrt{2|\Omega| \max\left((K_0)^2, (K_g)^2\right)}.$$

The variance bound $\mathrm{Var}(Z) \le C_Z$ holds uniformly for any $\alpha \in (0, 2]$ and $\beta \in (0, 1]$ because the boundedness of $u_0$, $g$, and the indicator functions ensures that the second moment is controlled, regardless of the specific distributions of $X^\alpha$ and $Y^\beta$. Consequently, $C$ is independent of $\alpha$ and $\beta$. □

**LEMMA 3.2.** *Let $B_2$ be the term in (3.9). Suppose that $f$ satisfies the condition (3.5). Consider the stochastic integral:*

$$\zeta := \int_0^{\tau_t \wedge \tau_\Omega(x)} f(t - Y^\beta(s), X^\alpha(s)) ds,$$

*and let $\{\zeta_j\}$ be i.i.d. copies of $\zeta$. Then the term $B_2$ in (3.9) is the Monte Carlo error:*

$$B_2 = \left\| \mathbb{E}[\zeta] - \frac{1}{M} \sum_{j=1}^{M} \zeta_j \right\|_{L^2(\Omega)},$$

*and it satisfies*

$$\left(\mathbb{E}[(B_2)^2]\right)^{1/2} \le \frac{C}{\sqrt{M}} \quad \text{with} \quad C = 2K_f \, diam(\Omega)^\alpha \sqrt{|\Omega|}.$$

*Note that the constant $C$ remains uniformly bounded for all $\alpha \in (0, 2]$ and $\beta \in (0, 1]$.*

*Proof.* Note that $\{\zeta_j\}$ denote i.i.d. copies of $\zeta$. Then the mean squared error is:

$$\mathbb{E}[(B_2)^2] \leq |\Omega| \sup_{x \in \Omega} \mathbb{E}\left[\left(\mathbb{E}[\zeta] - \frac{1}{M}\sum_{j=1}^M \zeta_j\right)^2\right] = |\Omega| \sup_{x \in \Omega} \frac{\text{Var}(\zeta)}{M}.$$

By the variance inequality $\text{Var}(\zeta) \leq \mathbb{E}[\zeta^2]$ and the the boundedness of $f$ in (3.5), we have

$$|\zeta| \leq \int_0^{\tau_t \wedge \tau_\Omega(x)} \left|f(t - Y^\beta(s), X^\alpha(s))\right| ds \leq K_f \cdot (\tau_t \wedge \tau_\Omega(x)),$$

which implies

$$\text{Var}(\zeta) \leq \mathbb{E}[\zeta^2] \leq (K_f)^2 \mathbb{E}\left[(\tau_t \wedge \tau_\Omega(x))^2\right] \leq (K_f)^2 \mathbb{E}\left[\tau_\Omega(x)^2\right].$$

According to (2.7) and Remark 2.1, we have

$$(3.12) \qquad\qquad \mathbb{E}\left[\tau_\Omega(x)^2\right] < \mathbb{E}\left[\tau_{\mathbb{B}_r^n}^2\right] \leq 4r^{2\alpha},$$

with $r$ being $\text{diam}(\Omega)$, the diameter of $\Omega$. As a result, we arrive at

$$\left(\mathbb{E}[(B_2)^2]\right)^{1/2} \leq \frac{C}{\sqrt{M}} \quad \text{with} \quad C = 2K_f \text{diam}(\Omega)^\alpha \sqrt{|\Omega|}$$

This completes the proof of the lemma. $\qquad\qquad\qquad\qquad\qquad\qquad\qquad\qquad\square$

**Lemma 3.3.** *Assume that $f$ satisfies the conditions (3.5) and (3.6). Then the term $B_3$ in (3.9) can be bounded by:*

$$(\mathbb{E}[(B_3)^2])^{1/2} \leq C(\Delta t)^{(1-\epsilon)/2},$$

*where the constant $C$ may depend on $t$, $\Omega$, $c_1$, $c_2$ and $\epsilon$, but is always independent of $\Delta t$, $\alpha$ and $\beta$.*

*Proof.* To begin, we define the error term

$$\Xi := \int_0^{\tau_t \wedge \tau_\Omega} f(s) ds - \Delta t \sum_{i=1}^k f(t_i) \quad \text{with } f(s) = f(t - Y^\beta(s), X^\alpha(s)),$$

and the corresponding pathwise error

$$\Xi_j := \int_0^{\tau_j} f_j(s) ds - \Delta t \sum_{i=1}^{k_j} f_j(t_i),$$

where $\tau_j = \tau_{t,j} \wedge \tau_{\Omega,j}(x)$ and $f_j(s) = f(t - Y_j^\beta(s), X_j^\alpha(s))$. Then we observe that

$$B_3 = \left\|\frac{1}{M}\sum_{j=1}^M \int_0^{\tau_j} f_j(s) ds - \frac{1}{M}\sum_{j=1}^M \left[\Delta t \sum_{i=1}^{k_j} f_j(t_i)\right]\right\|_{L^2(\Omega)} = \left\|\frac{1}{M}\sum_{j=1}^M \Xi_j\right\|_{L^2(\Omega)}.$$

Since the collection $\{(X_j^\alpha, \eta_j)\}_{j=1}^M$ consists of independent and identically distributed realizations of the random variable $\Xi$, we have

$$\mathbb{E}\left[\left(\frac{1}{M}\sum_{j=1}^M \Xi_j\right)^2\right] = \frac{\text{Var}(\Xi)}{M} + (\mathbb{E}[\Xi])^2 \leq \text{Var}(\Xi) + (\mathbb{E}[\Xi])^2 \leq \mathbb{E}[\Xi^2].$$

Assuming $\tau_t/\Delta t = N$ is an integer, then we observe that

$$\Xi = \sum_{i=1}^{k} \int_{t_{i-1}}^{t_i} [f(s) - f(t_i)] \, ds.$$

Now we consider the decomposition

$$f(s) - f(t_i) = \Big(f(s) - f(t - Y^\beta(t_i), X^\alpha(s))\Big) + \Big(f(t - Y^\beta(t_i), X^\alpha(s)) - f(t_i)\Big).$$

For the first term, we apply (3.6) to obtain

$$(3.13) \qquad |f(s) - f(t - Y^\beta(t_i), X^\alpha(s))| \le c_1 |Y^\beta(s) - Y^\beta(t_i)|^{\gamma_1}.$$

Meanwhile, for the second term, we apply (3.6) to obtain

$$(3.14) \qquad |f(t - Y^\beta(t_i), X^\alpha(s)) - f(t_i)| \le c_2 \|X^\alpha(s) - X^\alpha(t_i)\|^{\gamma_2}.$$

Thus

$$|\Xi| \le \sum_{i=1}^{k} \int_{t_{i-1}}^{t_i} c_1 |Y^\beta(s) - Y^\beta(t_i)|^{\gamma_1} + c_2 \|X^\alpha(s) - X^\alpha(t_i)\|^{\gamma_2} ds.$$

By the Young's inequality, we obtain

$$(3.15) \qquad \begin{aligned} \mathbb{E}[\Xi^2] &\le 2\mathbb{E}\left[\left(\sum_{i=1}^{k} \int_{t_{i-1}}^{t_i} c_1 |Y^\beta(s) - Y^\beta(t_i)|^{\gamma_1} ds\right)^2\right] \\ &\quad + 2\mathbb{E}\left[\left(\sum_{i=1}^{k} \int_{t_{i-1}}^{t_i} c_2 \|X^\alpha(s) - X^\alpha(t_i)\|^{\gamma_2} ds\right)^2\right]. \end{aligned}$$

For the first term, we derive

$$\begin{aligned} \mathbb{E}\left[\left(\sum_{i=1}^{k} \int_{t_{i-1}}^{t_i} c_1 |Y^\beta(s) - Y^\beta(t_i)|^{\gamma_1} ds\right)^2\right] &\le \mathbb{E}\left[(c_1)^2 k \sum_{i=1}^{k} \left(\int_{t_{i-1}}^{t_i} |Y^\beta(s) - Y^\beta(t_i)|^{\gamma_1} ds\right)^2\right] \\ &\le \mathbb{E}\left[(c_1)^2 t_k \sum_{i=1}^{k} \int_{t_{i-1}}^{t_i} |Y^\beta(s) - Y^\beta(t_i)|^{2\gamma_1} ds\right] \\ &\le \mathbb{E}\left[(c_1)^2 \tau_\Omega \sum_{i=1}^{k_\Omega} \int_{t_{i-1}}^{t_i} |Y^\beta(s) - Y^\beta(t_i)|^{2\gamma_1} ds\right] \\ &\le (c_1)^2 \mathbb{E}[(\tau_\Omega)^2] (\Delta t)^{-1} \int_0^{\Delta t} \mathbb{E}\left[|Y^\beta(s)|^{2\gamma_1}\right] ds, \end{aligned}$$

where in the last inequality, we use the fact that the random variables $\int_{t_{i-1}}^{t_i} |Y^\beta(s) - Y^\beta(t_i)|^{2\gamma_1} ds$ are i.i.d., and independent of $\tau_\Omega$. Now, according to (3.8), we use the fact that for $\gamma_1 < \beta/2 \le 1/2$

$$\mathbb{E}\left[|Y^\beta(s)|^{2\gamma_1}\right] \le \frac{4^{\gamma_2} \Gamma(\frac{1+2\gamma_1}{2}) \Gamma(1 - \frac{2\gamma_1}{\beta})}{\sqrt{\pi} \Gamma(1 - \gamma_1)} |s|^{2\gamma_1/\beta}.$$

Here we choose $\gamma_1 = \beta(1-\epsilon)/2$ with small $\epsilon > 0$ and obtain

$$\mathbb{E}\Big[|Y^\beta(s)|^{2\gamma_1}\Big] \leq C_\epsilon |s|^{1-\epsilon},$$

where $C_\epsilon$ is independent of $\beta$ but depends on $\epsilon$. This together with the estimate (3.12) leads to

$$(3.16) \qquad \mathbb{E}\Big[\Big(\sum_{i=1}^{k}\int_{t_{i-1}}^{t_i} c_1|Y^\beta(s) - Y^\beta(t_i)|^{\gamma_1}\mathrm{d}s\Big)^2\Big] \leq C_1\Delta t^{1-\epsilon}.$$

where the generic constant $C_1$ depends on $\epsilon$, $c_1$ and $\Omega$, but it is always independent of $\alpha$, $\beta$ and $\Delta t$. Similarly, we use the fact that $E[(\tau_t)^2]$ can be uniformly bounded with respect to $\beta$, and choose $\gamma_2 = \alpha(1-\epsilon)/2$ to obtain

$$(3.17) \qquad \mathbb{E}\left[\left(\sum_{i=1}^{k}\int_{t_{i-1}}^{t_i} c_2\|X^\alpha(s) - X^\alpha(t_i)\|^{\gamma_2}ds\right)^2\right] \leq C_2\Delta t^{1-\epsilon}.$$

where the constant $C_2$ depends on $\epsilon$, $c_2$ and $t$, but it is always independent of $\alpha$, $\beta$ and $\Delta t$.

Finally, Combining all bounds yields

$$\mathbb{E}\Big[(B_3)^2\Big] \leq (C_1 + C_2)|\Omega|(\Delta t)^{1-\epsilon}.$$

where $C_1$ and $C_2$ are constants appearing in (3.16) and (3.17). This completes the proof of the Lemma. □

Combining the estimates in Lemmas 3.1, 3.2 and 3.3 we conclude the following main theorem for the error estimates of the numerical scheme (3.3).

THEOREM 3.4. *Let $u$ and $u_M^*$ be the solutions of (1.1) and (3.3), respectively. Then there holds*

$$\Big(\mathbb{E}\Big[\|u(t,\cdot) - u_M^*(t,\cdot)\|_{L^2(\Omega)}^2\Big]\Big)^{\frac{1}{2}} \leq C\Big((\Delta t)^{(1-\epsilon)/2} + M^{-1/2}\Big).$$

*where $C$ is a positive constant independent of $\alpha$, $\beta$, $\Delta t$ and $M$.*

**4. Numerical results.** In this section, we provide several numerical examples to demonstrate the accuracy and robustness of the Monte Carlo method (3.3) in approximating the solution of the diffusion model (1.1). In the following, we apply the proposed algorithm to compute the numerical solution values $u_M^*(T,x)$ for any fixed $x \in \Omega$. To approximate the $L^2$-error, we randomly select a set of points $\{x_i\}_{i=1}^{N_s}$ and compute the error using the formula

$$\text{Error} = \|u(T,\cdot) - u_M^*(T,\cdot)\|_{L^2(\Omega)} \approx \sqrt{\frac{|\Omega|}{N_s}\sum_{i=1}^{N_s}\Big(u(T,x_i) - u_M^*(T,x_i)\Big)^2}.$$

Throughout our computations, we set $N_s = 1000$. A very fine time step is used to evaluate the convergence of the Monte Carlo approximation, while a sufficiently large $M$ is fixed to test the convergence of the time discretization.

**Example 1.** Let $\Omega = \mathbb{B}_1^n$ and $T = 1$. Consider the model (1.1) with the following exact solutions:

$$(4.1) \qquad u(t,x) = t^\beta(1 - |x|^2)_+^{\frac{\alpha}{2}},$$

where $b_+ = \max\{b, 0\}$. In this problem, $g(t, x) = 0$ on $\Omega^c$ and the source term is

$$f(t, x) = \Gamma(\beta + 1)(1 - |x|^2)_+^{\frac{\alpha}{2}} + t^\beta 2^\alpha \frac{\Gamma(1 + \frac{\alpha}{2})\Gamma(\frac{n+\alpha}{2})}{\Gamma(\frac{n}{2})}.$$

Note that both the exact solution and the source data exhibit a weak singularity near $t = 0$ and a weakly singular boundary layer near $\partial\Omega$. These singularities become more pronounced as $\alpha \to 0^+$ and $\beta \to 0^+$.

First, we set $n = 2$ and test the error. As we mentioned earlier, to evaluate the convergence of the Monte Carlo approximation, we fix a very fine time step of $\Delta t = 10^{-4}$. In Figure 2, we plot the $L^2$ errors of the proposed stochastic method for different numbers of paths $M$. As shown in Figure 2, the error decreases as the number of paths $M$ increases, and the convergence rate aligns with the estimates in Theorem 3.4, approximately of order $O(M^{-1/2})$.

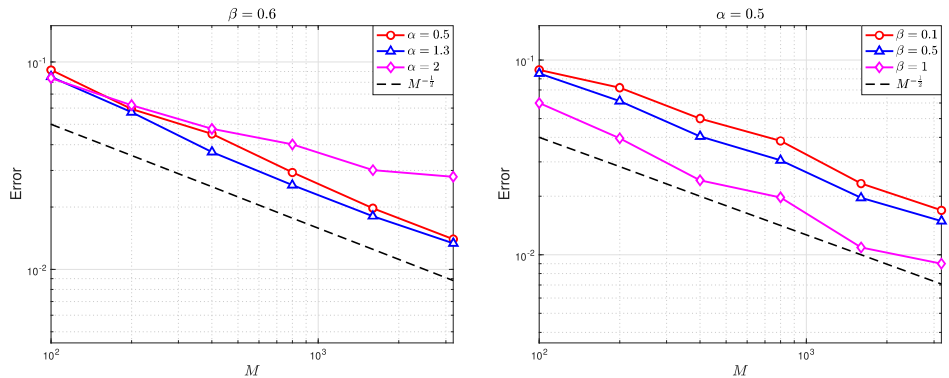

Fig. 2. *Example 1: errors against $M$ (the number of paths) with various $\alpha$ (left) and $\beta$ (right).*

Next, in Figure 3, we fix $M = 10^4$ and examine the convergence of the time discretization. The left panel of Figure 3 shows the $L^2$-error plotted against various time step sizes $\Delta t$ for $\alpha = 0.5$, 1.3, and 2.0, with $\beta$ fixed at 0.6. The right panel presents the $L^2$-error for different time step sizes $\Delta t$, considering $\beta = 0.1$, 0.5, and 1.0, while keeping $\alpha = 0.5$. The numerically observed convergence rates of the $L^2$-error align with the error estimates established in Theorem 3.4, demonstrating an approximate half-order rate.

One of the key advantages of the proposed stochastic algorithm is its feasibility in high-dimensional settings. To evaluate this, we test cases with $n = 100$, which are particularly challenging for traditional numerical methods. Additionally, we consider another demanding scenario where $\beta$ is very close to zero, specifically $\beta = 0.03$. In Figure 4, we report the error for different values of $T$. The results confirm convergence rates of $O(M^{-\frac{1}{2}})$ and $O(\Delta t^{\frac{1}{2}})$, as predicted by Theorem 3.4. These findings underscore the effectiveness and robustness of the proposed stochastic algorithm, even in scenarios with small $\beta$ and high-dimensional problems.

To illustrate the robustness of our method with respect to $\alpha$, we test the error convergence rate for different $T$ with $\alpha = 0.02$ and $n = 100$. In the left panel of Figure 5, the error decreases as the number of paths $M$ increases, displaying a convergence rate of approximately $M^{-1/2}$ for a fixed time step $\Delta t = 10^{-3}$. Similarly, in the right panel, with $M = 10^4$, the error decreases at a stable rate sligtly higher than $O(\Delta t^{1/2})$. These results are fully support our theory in Theorem 3.4, showcasing the effectiveness and robustness of our method for very small $\alpha$, even in high-dimensional settings. Additionally, Figure 6 presents numerical results for the case where both $\alpha$ and $\beta$ diminish simultaneously, specifically $\alpha = 0.05$ and $\beta = 0.05$. The convergence remains robust in this challenging scenario, highlighting one of the key advantages of the proposed method. This robustness distinguishes our approach from other existing methods for solving fractional PDEs.

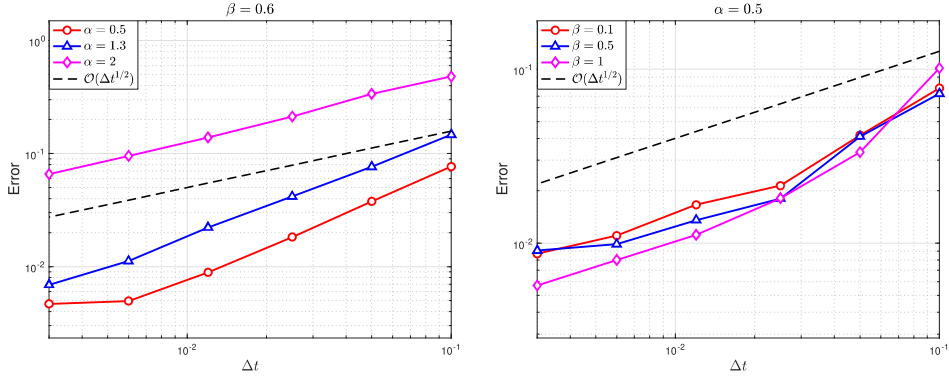

Fig. 3. *Example 1: errors against time step $\Delta t$ with various $\alpha$ (left) and $\beta$ (right).*

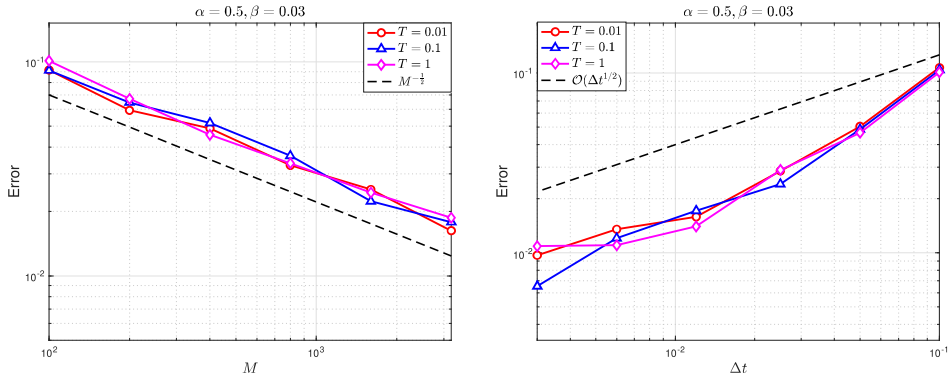

Fig. 4. *Example 1 (n=100): errors against $M$ (left) and $\Delta t$ (right), with $\alpha = 0.5$ and $\beta = 0.03$.*

**Example 2.** Let $\Omega = \mathbb{B}_1^2$ and $T = 1$. We consider the model (1.1) with the following exact solution:

$$(4.2) \qquad u(t,x) = E_{\beta,1}(-t^\beta)(1 - |x|^2)_+^{\frac{\alpha}{2}},$$

where $E_{\beta,1}$ represents the Mittag–Leffler function. In this case, the initial condition is $u_0 = (1 - |x|^2)_+^{\frac{\alpha}{2}}$, the source term is

$$f(t,x) = -E_{\beta,1}(-t^\beta)(1 - |x|^2)_+^{\frac{\alpha}{2}} + E_{\beta,1}(-t^\beta)2^\alpha \Gamma^2(1 + \alpha/2),$$

and the boundary condition is the homogeneous condition $g \equiv 0$. Compared with **Example 1**, the function (4.2) exhibits more complex singularities in the time variable as $t \to 0^+$.

First, we fix the time step to $\Delta t = 5 \times 10^{-4}$. The $L^2$ error of the proposed scheme for different numbers of paths $M$ is shown in Figure 7. We observe that the error decreases as the number of paths $M$ increases, with a convergence rate approximately of the order $M^{-1/2}$. Next, we test the convergence rates of the time discretization for the proposed algorithm by fixing $M = 10^4$. In the left panel of Figure 8, we present the $L^2$ error of the numerical solution for various time step sizes $\Delta t$ with $\alpha = 0.5$, 1, 2.0 and $\beta = 0.6$. In the right panel of Figure 8, we display the error of the numerical solution (3.3) for various time step sizes $\Delta t$ with $\beta = 0.1$, 0.5, 1 and $\alpha = 0.5$. It is noted that the numerical errors decrease as the time step $\Delta t$ decreases, and the convergence rate of the errors under the $L^2$ norm is approximately of half-order for different values of $\alpha$ and $\beta$, which is consistent with the estimates in Theorem 3.4.

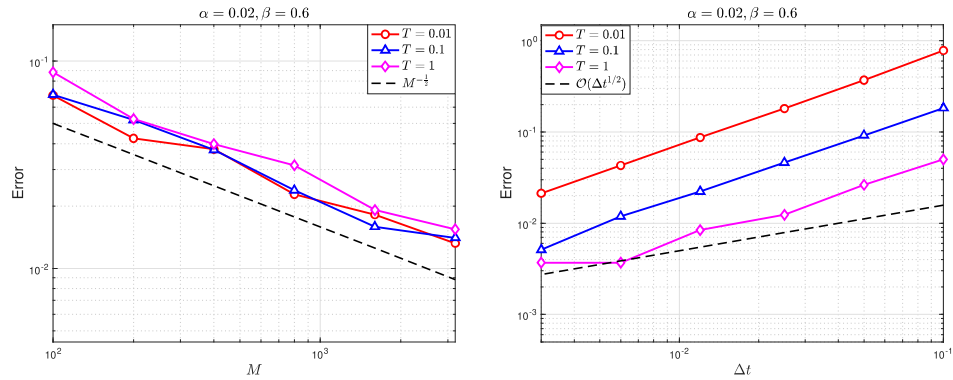

FIG. 5. *Example 1 (n = 100): errors against M (left) and Δt (right), with α = 0.02 and β = 0.6.*

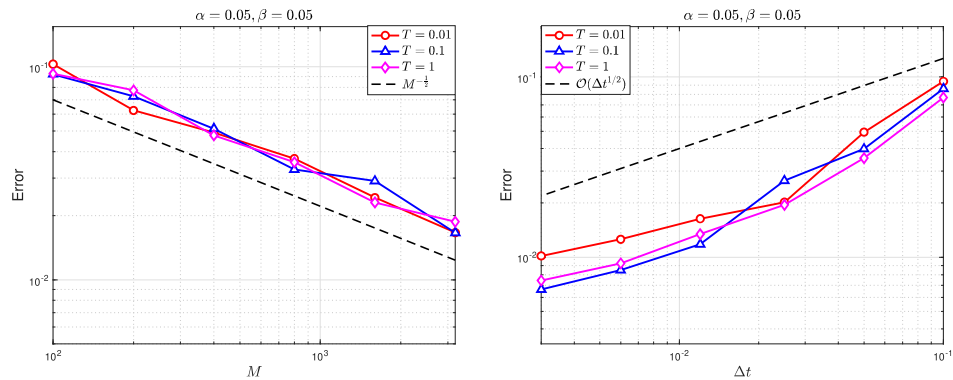

FIG. 6. *Example 1 (n = 100): errors against M (left) and Δt (right), with α = 0.05 and β = 0.05.*

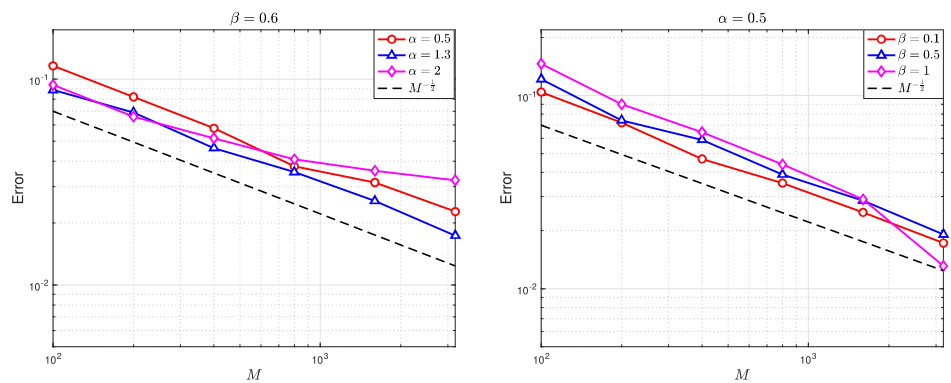

FIG. 7. *Example 2: errors against M (the number of paths) with various α (left) and β (right).*

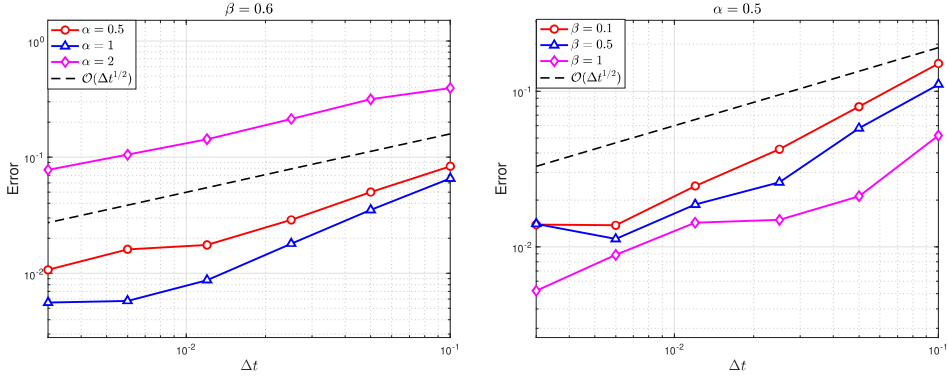

Fig. 8. *Example 2: errors against time step $\Delta t$ with various $\alpha$ (left) and $\beta$ (right).*

**Example 3.** In the previous examples, the boundary condition was set to zero. Here, we present numerical experiments with a nonhomogeneous boundary condition, i.e., $g \neq 0$. Specifically, we consider the two-dimensional problem (1.1) on an L-shaped domain $\Omega = (-1, 1)^2 \setminus [0, 1]^2$, with the exact solution defined as

$$u(t, x) = t^a (1 + |x|^2)^{-\frac{7}{2}} \quad \text{with } a = 1.2.$$

The boundary condition is specified as $g(t, x) = u(t, x)$ for all $(t, x) \in (0, T] \times \Omega^c$. This condition is nonzero and varies with time. The source term is

$$f(t, x) = \frac{\Gamma(a+1) t^{a-\beta}}{\Gamma(a+1-\beta)} (1 + |x|^2)^{-\frac{7}{2}} + \frac{2^\alpha \Gamma\left(\frac{\alpha+7}{2}\right) \Gamma\left(\frac{\alpha+2}{2}\right) t^a}{\Gamma\left(\frac{7}{2}\right)} {}_2F_1\left(\frac{2+\alpha}{2}, \frac{7+\alpha}{2}; 1; -|x|^2\right).$$

To test the convergence with respect to the number of paths $M$, we fix the time step as $\Delta t = 5 \times 10^{-4}$. The results in Figure 9 illustrate the errors of the numerical scheme for various values of $M$. As observed, the error decreases with $M$ at a rate of $O(M^{-1/2})$, consistent with the estimates in Theorem 3.4. Similarly, we test the convergence rate of the numerical solution with respect to $\Delta t$, fixing the number of paths at $M = 10^4$. Figure 10 shows that the convergence rate is $O(\Delta t^{1/2})$, which confirms our theoretical result.

**Example 4.** We now consider the hexagonal hailstone domain to further show the efficiency and robustness of our method. In the computation, we adopt the source term

$$f(t, x) = \frac{\cos(t)}{1 + 10t^2} \left( \cos(\frac{\pi}{3} x_1^2 - x_1 x_2) + \sin(\frac{\pi}{6} x_2^2 + x_1 x_2) \right),$$

and the boundary condition $g \equiv 0$. The corresponding computational domain $\Omega$ can be determined by the following polar coordinate transformation of the form

$$\begin{cases} x = R(\theta) \cos(\theta) \\ y = R(\theta) \sin(\theta). \end{cases}$$

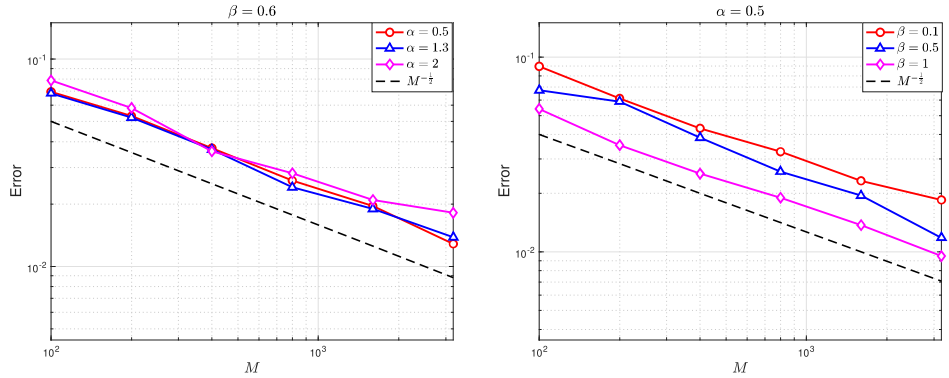

Fig. 9. *Example 3: errors against M (the number of paths) with various α (left) and β (right).*

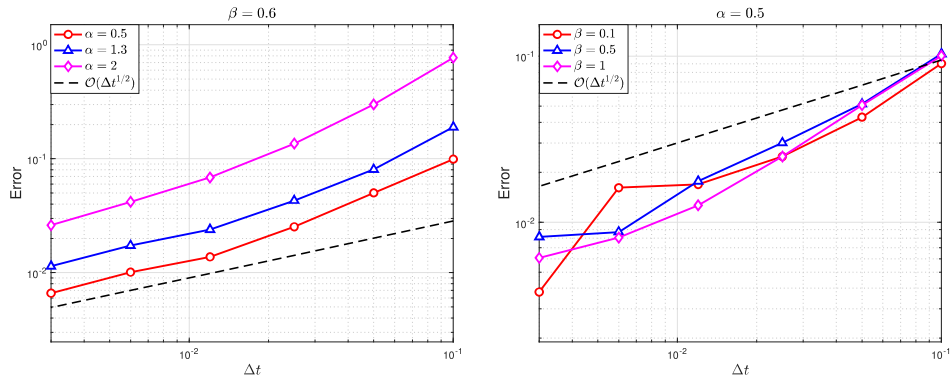

Fig. 10. *Example 3: errors against time step Δt with various α (left) and β (right).*

Here, we take $R(\theta) = 1 + 0.9 \sin(6\theta) + 0.1 \cos(10\theta)$ with $\theta \in [0, 2\pi]$, and the initial value is chosen as the uniformly distributed data between the range $[-2, -1]$.

We present the profiles of the numerical solutions for various fractional orders $\alpha$ in Figure 11, with $T = 0.5$ and $\beta = 0.6$ fixed. As expected, the singularity layer near the boundary becomes thinner as the fractional power $\alpha$ decreases. Next, Figure 12 shows the profiles of the numerical solutions for different values of $\beta$, with $T = 1$ and $\alpha = 1$ fixed. As $\beta$ decreases, the singularity layer near the boundary becomes thicker. This increased non-locality induces long-memory behavior in the solution, leading to an expansion of the boundary layer. Furthermore, Figure 13 illustrates the evolution of the numerical solution over time, with $\alpha = 1$ and $\beta = 0.5$ fixed. The boundary layer gradually thickens over time and eventually stabilizes at a steady state. This example demonstrates the feasibility and robustness of the proposed stochastic scheme, even for domains with complex geometries.

## 5. Concluding remarks.
Space-time fractional PDEs in high-dimensional space have found wide applications in science and engineering. The numerical method based on the Feynman–Kac formula is very popular and simple to solve this kind of high-dimensional problem. This paper proposed a Monte Carlo method based on the Feynman–Kac formula for solving space-time fractional PDEs on complex domains in multiple dimensions. The key idea is to combine the simulation of the monotone path of a stable subordinator in time corresponding to the time-fractional direvative with the "walk-on-spheres" method that efficiently simulates the stable Lévy jumping process in space associate with the fractional Laplacian. We derived the error estimate of the proposed scheme, which is affected by the time step size $O(\Delta t^{(1-\epsilon)/2})$ and the number of paths $O(M^{-\frac{1}{2}})$. The numerical experiments confirmed the theoretical expectations and demonstrated that the suggested algorithm allows us to use parallel computing to improve the calculation efficiency and significantly reduce the computation time. Furthermore, the robustness of the scheme across various fractional orders highlights its potential for simulating local-nonlocal and nonlocal-nonlocal coupling models, which remain an open area for future research.

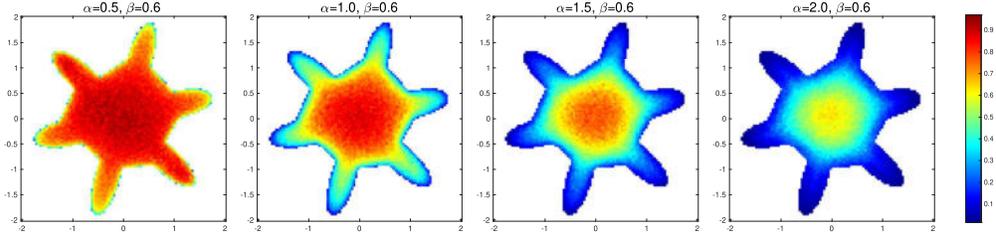

Fig. 11. *Example 4: Numerical solutions with $\beta = 0.6$ and various $\alpha$ at $T = 0.5$.*

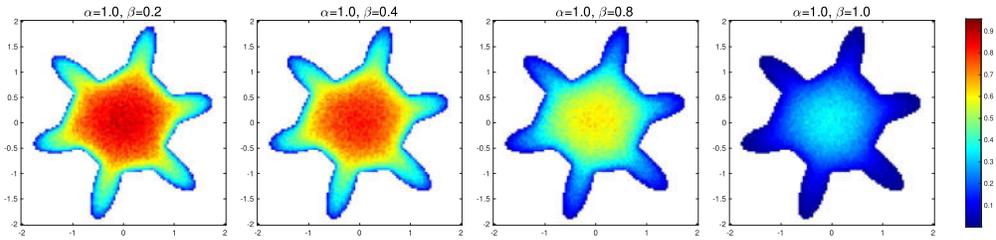

Fig. 12. *Example 4: Numerical solutions with $\alpha = 1$ and various $\beta$ at $T = 1$.*

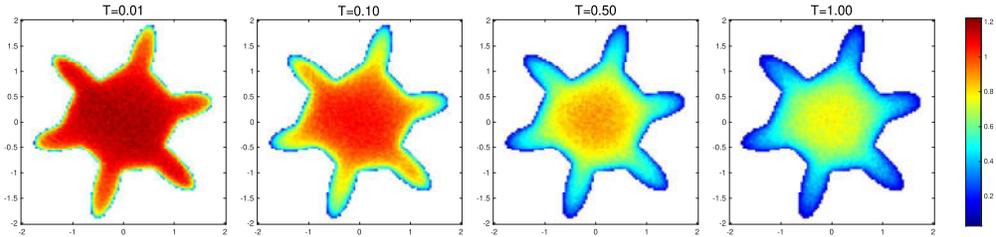

Fig. 13. *Example 4: Numerical solutions at various $T$ with fixed fractional order $\alpha = 1$ and $\beta = 0.5$.*